\documentclass{article}%
\usepackage{amsfonts}
\usepackage{amsmath}
\usepackage{amssymb}
\usepackage{graphicx}%
\providecommand{\U}[1]{\protect \rule{.1in}{.1in}}
\newtheorem{theorem}{Theorem}

\newtheorem{conjecture}[theorem]{Conjecture}
\newtheorem{corollary}[theorem]{Corollary}

\newtheorem{definition}[theorem]{Definition}

\newtheorem{proposition}[theorem]{Proposition}
\newtheorem{remark}[theorem]{Remark}

\newenvironment{proof}[1][Proof]{\noindent \textbf{#1.} }{\  \rule{0.5em}{0.5em}}
\begin{document}

\begin{center}
\textbf{Linear maps preserve the covariance sets }

\bigskip

{\large Mohammad Hossein Alizadeh}\footnote{e-mail: alizadeh@aegean.gr}

\bigskip

{\large Department of Mathematics,}

{\large Islamic Azad University, Nur branch, Nur, Iran}

\bigskip
\end{center}

\textbf{Abstract. }Let $\mathcal{A}$ and $\mathcal{B}$\ are $C^{{\huge \ast}}%
$-algebras\textbf{.} A\textbf{ }linear map $\phi:\mathcal{A\rightarrow B}$ is
$C^{\ast}$-Jordan homomorphism if it is a Jordan homomorphism which preserves
the adjoint operation. In this note we show that $C^{\ast}$-Jordan
homomorphisms -under mild assumptions- preserving covariance set and
covariance coset in $C^{\ast}$-algebras.

\bigskip


\bigskip


\vfill

\pagebreak

\section{Introduction and preliminaries}

Suppose that $\mathcal{A}$ is a $C^{\ast}$-algebra with identity $1$. An
element $a\in \mathcal{A}$ is called \textit{regular} if it has a generalized
inverse in $\mathcal{A}$, i.e. there exists $b\in \mathcal{A}$ such that%

\[
aba=a.
\]

We say that an element $a\in \mathcal{A}$\textit{\ is Moore-Penrose invertible
if there exists }$b\in \mathcal{A}$\textit{\ such that}%
\[
aba=a,\  \  \ bab=b,\  \  \  \left(  ab\right)  ^{\ast}=ab\  \  \  \text{and
\  \ }\left(  ba\right)  ^{\ast}=ba.
\]

It is well known that the Moore-Penrose inverse (briefly, MP--inverse) is
unique if it exists. We reserve the notation $a^{\dag}$ for the MP--inverse of
$a.$ In what follows, we will denote by $\mathcal{A}^{-1}$ the subset of
invertible elements of $\mathcal{A}$ and by $\mathcal{A}^{\dag}$, the set of
all MP--invertible elements of $\mathcal{A}$. The \textit{commutator} of a
pair of elements $x$ and $y$ in $\mathcal{A}$\ is given by
\[
\left[  x,y\right]  =xy-yx.
\]
Note that $\left[  x,y\right]  =0$ if and only if $x$ and $y$ commute.

In the next section we need the following definition of covariance set which
was studied in \cite{Aliz2011}

\begin{definition}
\cite{Aliz2011} For a given element $a\in \mathcal{A}^{\dag}$ with MP--inverse
$a^{\dag}$\ we will denote the \textit{covariance }set by $\mathfrak{C}(a)$
and define;
\begin{equation}
\mathfrak{C}(a)=\{b\in \mathcal{A}^{-1}:(bab^{-1})^{\dag}=ba^{\dag}b^{-1}\}
\text{.} \label{coset}%
\end{equation}

\end{definition}

Also the notion of covariance coset was introduced and studied in
\cite{Aliz2013}. In fact, this set is defined by reversing the roles of $a$
and $b$ in $\mathfrak{C}(a)$ and is denoted by $\mathfrak{B}\left(  b\right)
$. i.e.,%

\begin{equation}
\mathfrak{B}\left(  b\right)  =\left \{  a\in \mathfrak{A}^{\dag}:(bab^{-1}%
)^{\dag}=ba^{\dag}b^{-1}\right \}  . \label{coco}%
\end{equation}

The propose of this work is to show that under weak assumptions, $C^{\ast}%
$-Jordan homomorphisms preserving covariance set and covariance coset in
$C^{\ast}$-algebras.

\section{Main results}

We recall the following definitions and theorems which will be needed to prove
some of our results.

\begin{definition}
\cite{BM2010} We say that a $C^{\ast}$-algebra $\mathcal{A}$ is of real rank
zero if the set formed by all the real linear combinations of (orthogonal)
projections is dense in the set of self-adjoint elements of $\mathcal{A}$.
\end{definition}

\begin{remark}
\label{Rem1}Suppose that $\mathcal{A}$ and $\mathcal{B}$ are $C^{\ast}%
$-algebras. It is well known that (see \cite{BM2010}) the property of the
above definition is satisfied by every von Neumann algebra, and in particular
by the $C^{\ast}$-algebra $B(H)$ of all bounded linear operators on a Hilbert
space $H$, and by the Calkin algebra $C(H)=\frac{B(H)}{K(H)}$.
\end{remark}

\begin{definition}
We say that a linear map $\phi:\mathcal{A\rightarrow B}$ is $C^{\ast}$-Jordan
homomorphism if it is a, Jordan homomorphism which preserves the adjoint
operation, i.e.%
\[
\phi \left(  x^{\ast}\right)  =\left(  \phi \left(  x\right)  \right)  ^{\ast
}\qquad \forall x\in \mathcal{A}\text{.}%
\]

\end{definition}

The $C^{\ast}$-homomorphism and $C^{\ast}$-anti-homomorphism are analogously defined.

In 2012, Boudi and Mbekhta \cite{BM2010}\ proved the following theorem.

\begin{theorem}
\label{Mbekh}Let $\mathcal{A}$ be a $C^{\ast}$-algebra of real rank zero and
$\mathcal{B}$ a prime $C^{\ast}$ -algebra. Let $\phi:\mathcal{A\rightarrow B}$
be a surjective, unital linear map. Then the following conditions are equivalent:

1) $\phi(x^{\dag})=\left(  \phi(x)\right)  ^{\dag}$ for all $x\in A^{\dag};$

2) $\phi$ is either a $C^{\ast}$-homomorphism or a $C^{\ast}$-anti-homomorphism.
\end{theorem}

\begin{proof}
See \cite[Theorem 3.3]{BM2010}.
\end{proof}

The next proposition describes a relation between the covariance set
$\mathfrak{C}(a),$ and commutators. It was proved in \cite{Aliz2011}.

\begin{proposition}
\label{comm}Let $a$ $\in$ $\mathcal{A}^{\dag}$ with MP--inverse $a^{\dag}.$
Then the following statements are equivalent:

(i) $b\in \mathfrak{C}(a);$

(ii) \ $[b^{\ast}b,aa^{\dag}]=0$ and $[b^{\ast}b,a^{\dag}a]=0$.
\end{proposition}

A similar result also is true for covariance coset:

\begin{proposition}
\cite{Aliz2013} \label{WMP-equv}Assume $b\in \mathfrak{A}^{-1}$. Then the
following statements are equivalent:

(i)\ $a\in \mathfrak{B}\left(  b\right)  ;$

(ii) $\left[  a^{\dag}a,b^{\ast}b\right]  =0$ \textit{and }$\left[  aa^{\dag
},b^{\ast}b\right]  =0.$
\end{proposition}

Now we are going to prove the main result.

\begin{theorem}
\label{cofi}Let $\mathcal{A}$ be a $C^{\ast}$-algebra of real rank zero and
$\mathcal{B}$ a prime $C^{\ast}$ -algebra. Let $\phi:\mathcal{A\rightarrow B}$
be a surjective, unital linear map. If $\phi(x^{\dag})=\left(  \phi(x)\right)
^{\dag}$ for all $x\in A^{\dag},$ then $\phi \left(  \mathfrak{C}(a)\right)  =$
$\mathfrak{C}(\phi \left(  a\right)  )$ and $\phi \left(  \mathfrak{B}%
(a)\right)  =$ $\mathfrak{B}(\phi \left(  a\right)  ).$
\end{theorem}

\begin{proof}
By Theorem \ref{Mbekh}, $\phi$ is either a $C^{\ast}$-homomorphism or a
$C^{\ast}$-anti-homomorphism. First we assume that $\phi$ is a $C^{\ast}%
$-homomorphism. Let $b\in \mathfrak{C}(a)$. By Proposition \ref{comm}
\begin{equation}
b^{\ast}baa^{\dag}=aa^{\dag}b^{\ast}b,\qquad b^{\ast}ba^{\dag}a=a^{\dag
}ab^{\ast}b \label{1}%
\end{equation}
Since $\phi$ is a $C^{\ast}$-homomorphism, from (\ref{1}) we get%
\begin{align*}
\phi \left(  b\right)  ^{\ast}\phi \left(  b\right)  \phi \left(  a\right)
\phi \left(  a\right)  ^{\dag}  &  =\phi \left(  a\right)  \phi \left(  a\right)
^{\dag}\phi \left(  b\right)  ^{\ast}\phi \left(  b\right)  ,\\
\phi \left(  b\right)  ^{\ast}\phi \left(  b\right)  \phi \left(  a\right)
^{\dag}\phi \left(  a\right)   &  =\phi \left(  a\right)  ^{\dag}\phi \left(
a\right)  \phi \left(  b\right)  ^{\ast}\phi \left(  b\right)
\end{align*}
which means that $\phi \left(  b\right)  \in \phi \left(  \mathfrak{C}(a)\right)
$ i.e. $\phi \left(  \mathfrak{C}(a)\right)  \subset \mathfrak{C}(\phi \left(
a\right)  )$. Since $\phi$ is surjective we get $\phi \left(  \mathfrak{C}%
(a)\right)  =\mathfrak{C}(\phi \left(  a\right)  ).$

Now we suppose that $\phi$ is a $C^{\ast}$-anti-homomorphism. Let
$b\in \mathfrak{C}(a)$. Again by Proposition \ref{comm} we have (\ref{1}).
Applying $\phi$ on (\ref{1}) we get
\begin{equation}
\phi \left(  aa^{\dag}\right)  \phi \left(  b^{\ast}b\right)  =\phi \left(
b^{\ast}b\right)  \phi \left(  aa^{\dag}\right)  ,\quad \phi \left(  a^{\dag
}a\right)  \phi \left(  b^{\ast}b\right)  =\phi \left(  b^{\ast}b\right)
\phi \left(  a^{\dag}a\right)  . \label{2}%
\end{equation}
Since $\phi$ is a $C^{\ast}$-anti-homomorphism and $\phi(x^{\dag})=\left(
\phi(x)\right)  ^{\dag}$ from (\ref{2}) we obtain%
\begin{align*}
\phi \left(  a)^{\dag}\phi(a\right)  \phi \left(  b)\phi(b^{\ast}\right)   &
=\phi \left(  b)\phi(b^{\ast}\right)  \phi \left(  a)^{\dag}\phi(a\right)  ,\\
\phi \left(  a\right)  \phi(a)^{\dag}\phi \left(  b)\phi(b^{\ast}\right)   &
=\phi \left(  b)\phi(b^{\ast}\right)  \phi \left(  a\right)  \phi(a)^{\dag}%
\end{align*}
Now by using Proposition \ref{comm} we conclude that $\phi \left(  b\right)
\in \phi \left(  \mathfrak{C}(a)\right)  $ i.e. $\phi \left(  \mathfrak{C}%
(a)\right)  =\mathfrak{C}(\phi \left(  a\right)  )$.

Applying Proposition \ref{WMP-equv}, a similar argument shows that
$\phi \left(  \mathfrak{B}(a)\right)  =\mathfrak{B}(\phi \left(  a\right)  )$.
\end{proof}

By Theorem \ref{cofi} and Remark \ref{Rem1}, we deduce the following results.

\begin{corollary}
Assume that $\mathcal{A}$ and $\mathcal{B}$ are $C^{\ast}$-algebras and also
von Neumann algebras. Let $\phi:\mathcal{A\rightarrow B}$ be a surjective,
unital linear map. If $\phi(x^{\dag})=\left(  \phi(x)\right)  ^{\dag}$ for all
$x\in A^{\dag},$ then $\phi \left(  \mathfrak{C}(a)\right)  =$ $\mathfrak{C}%
(\phi \left(  a\right)  )$ and $\phi \left(  \mathfrak{B}(a)\right)  =$
$\mathfrak{B}(\phi \left(  a\right)  ).$
\end{corollary}

\begin{corollary}
Suppose that $H$ and $K$ are Hilbert spaces. Let $\phi:B\left(  H\right)
\mathcal{\rightarrow}B\left(  K\right)  $ be a surjective linear map. If
$\phi(T^{\dag})=\left(  \phi(T)\right)  ^{\dag}$ for all $T\in B\left(
H\right)  ^{\dag},$ then $\phi \left(  \mathfrak{C}(T)\right)  =$
$\mathfrak{C}(\phi \left(  T\right)  )$ and $\phi \left(  \mathfrak{B}%
(T)\right)  =$ $\mathfrak{B}(\phi \left(  T\right)  ).$
\end{corollary}

Let $n\in%
\mathbb{N}
$. We say that a linear map $\phi:\mathcal{A\rightarrow B}$ is $n$-$C^{\ast}%
$-Jordan homomorphism if it is a, $n$-Jordan homomorphism (for more detail see
\cite{Eshaghi 2009}) which preserves the adjoint operation.

\textbf{Question}: For which $n\in%
\mathbb{N}
,$ the above results are true for $n$-$C^{\ast}$-Jordan homomorphism?

In connection with Theorem \ref{cofi}, we conclude the paper by the following conjecture:

\begin{conjecture}
Assume that $\mathcal{A}$ and $\mathcal{B}$ are $C^{\ast}$-algebras. Let
$\phi:\mathcal{A\rightarrow B}$ be a surjective, unital linear map. If
$\phi(x^{\dag})=\left(  \phi(x)\right)  ^{\dag}$ for all $x\in A^{\dag},$ then
$\phi \left(  \mathfrak{C}(a)\right)  =$ $\mathfrak{C}(\phi \left(  a\right)  )$
and $\phi \left(  \mathfrak{B}(a)\right)  =$ $\mathfrak{B}(\phi \left(
a\right)  ).$
\end{conjecture}

\end{document}